\UseRawInputEncoding
\documentclass{amsart}
\usepackage{amsmath}
\usepackage{amssymb}
\newtheorem{thm}{Theorem}[section]
\newtheorem{cor}[thm]{Corollary}
\newtheorem{conj}[thm]{Conjecture}
\newtheorem{lem}[thm]{Lemma}
\newtheorem{prop}[thm]{Proposition}
\newtheorem{cons}[thm]{Construction}
\theoremstyle{definition}
\newtheorem{defn}[thm]{Definition}
\theoremstyle{remark}
\newtheorem{rem}[thm]{Remark}
\newtheorem{ex}[thm]{Example}
\long\def\Thm#1{\begin{thm} #1 \end{thm}}
\long\def\Cor#1{\begin{cor} #1 \end{cor}}

\long\def\Lem#1{\begin{lem} #1 \end{lem}}
\long\def\Prop#1{\begin{prop} #1 \end{prop}}

\long\def\Rem#1{\begin{rem} #1 \end{rem}}

\def\bar#1{\overline{#1}}
\def\Sect{\section}
\def\Rarr#1#2{\xrightarrow[#2]{#1}}

\long\def\Ref#1#2#3#4#5#6{
\bibitem{#1}
{\rm #2,}
\textit{#3.}
{\rm #4}
\textbf{#5}
{\rm #6.}
}
\long\def\Refb#1#2#3#4{
\bibitem{#1}
{\rm #2,}
\textit{#3.}
#4.
}
\def\Zz{{\mathbb Z}}
\def\Rr{{\mathbb R}}
\def\Cc{{\mathbb C}}
\def\Tt{{\mathbb T}}
\def\Ff{{\mathbb F}}

\def\Sym{{\rm S}}
\def\i{{\rm i}}

\def\e{{\rm e}}
\def\phi{\varphi}
\def\into{\hookrightarrow}
\def\leq{\leqslant}
\def\geq{\geqslant}
\def\comp{\mathbin{\mathchoice
{\circ}
{{\scriptstyle\circ}}
{{\scriptscriptstyle\circ}}
{{\scriptscriptstyle\circ}}
}}
\def\st{\mid}
\def\Zero{{\rm Zero}}
\def\Null{{\rm Null}}

\def\semidirect{\rtimes}
\def\map{{\rm map}}
\def\cl#1{\overline{#1}}
\begin{document}

\title{A Borsuk--Ulam theorem for cyclic $p$-groups}
\author{M.~C.~Crabb}
\address{%
Institute of Mathematics\\
University of Aberdeen \\
Aberdeen AB24 3UE \\
UK}
\email{m.crabb@abdn.ac.uk}
\date{November 2022}
\begin{abstract}
We describe a connective $K$-theory Borsuk--Ulam/Bourgin--Yang theorem
for cyclic groups of order a power of a prime $p$.
Consider two finite dimensional complex representations $U$ and $V$
of the cyclic group $\Zz /p^{k+1}$ of order $p^{k+1}$, where $k\geq 0$.
For $0\leq l\leq k$, we write $V_l$ for the subspace of $V$ fixed by
the cyclic subgroup of order $p^l$, and
require that the fixed subspace, $V_{k+1}$, be zero 
and that $V_k$ be non-zero. Put
$\delta (V)=\sum_{l=0}^k p^l \dim _\Cc (V_l/V_{l+1})-(p^k-1)$.
Then the zero-set of any $\Zz /p^{k+1}$-map
$S(U) \to V$ from the unit sphere in $U$ (for some invariant
inner product) has covering dimension greater than or equal to
$2(\dim_\Cc U - \delta (V)-1)$, if $\dim_\Cc U> \delta (V)$.
\end{abstract}
\maketitle
\Sect{Introduction}
Fixing a prime $p$ and a natural number $k\geq 0$, we consider
the cyclic group $\Zz /p^{k+1}$ of order $p^{k+1}$.

The symbol $V$ will be reserved for
an $m$-dimensional complex representation
of $\Zz /p^{k+1}$ with trivial fixed submodule
$V^{\Zz /p^{k+1}}=0$, and,
for $0\leq l\leq k$, we set 
$$
m_l=\dim_\Cc (V^{p^{k-l+1}\Zz /p^{k+1}}/V^{p^{k-l}\Zz /p^{k+1}}),
$$
so that $m=\sum_{l=0}^k m_l$.
If $m_k\not=0$, that is, $V^{p\Zz /p^{k+1}}\not=0$, we define
$$
\delta (V) = \sum_{l=0}^k p^lm_l -(p^k-1).
$$
Notice that $\delta (V)\geq \dim V$, with equality only if 
either $k=0$ (and so $m_0=m$) 
or $k\geq 1$ and $m_0=m-1$ (and so $m_k=1$).

In the statement of the first theorem, and throughout the paper, we use
the representable cohomology, $H^*$, for compact Hausdorff spaces
(as, for example, in \cite[Section 8]{borsuk}).
\Thm{\label{main}
Let $X$ be a compact connected free $\Zz /p^{k+1}$-ENR such that, 
for some integer $n>1$,
the $\Ff_p$-homology group
$H_i(X;\,\Ff_p)$ is zero in dimensions $0<i<2(n-1)$.
Suppose that $m_k\not=0$ and $n>\delta (V)$.
Consider a $\Zz /p^{k+1}$-equivariant map
$f : X \to V$.

Then the zero-set
$$
\Zero (f) = \{ x \in X \st f(x)=0\}
$$
is non-empty and there exists
an integer $d\geq 2(n-1-\delta (V))$ such that 
the mod $p$ cohomology group $H^d(\Zero (f);\, \Ff_p)$ 
is non-zero; hence the covering dimension of $\Zero (f)$ is 
greater than or equal to $2(n-1-\delta (V))$.
}
\Rem{\label{red}
If for some $k'$, $0\leq k'<k$, $m_l=0$ for $l>k'$, and $m_{k'}\not=0$, 
we can restrict from $\Zz /p^{k+1}$ to the
subgroup $p^{k-k'}\Zz /p^{k+1}$ of order $p^{k'+1}$
and so, defining $\delta (V)$ in this case to be
$$
\delta (V)=\sum_{l=0}^{k'}p^lm_l-(p^{k'}-1),
$$
we again have the bound $2(n-1-\delta (V))$.
}
Results of this type go back to the work of Conner and Floyd 
\cite[(33.1)]{cf}
and Munkholm \cite{munkholm0, munkholm} some fifty years 
ago, but have attracted continuing attention (notably
in \cite{bartsch, ko,
MMS, BMS}) over the intervening years.

The bound on the dimension in Theorem \ref{main}
is close to being the best possible.
Let $L$ denote the basic $1$-dimensional complex representation
of $\Zz /p^{k+1}$:
$\Cc$ with the generator $1$ acting as multiplication by 
$\zeta=\e^{2\pi\i/p^{k+1}}$.
For $n\geq 1$, we write $nL=L\otimes_\Cc \Cc^n$
and $S(nL)$ for the unit sphere in $nL$.

It follows from Theorem \ref{main}
(and was already shown in \cite{ko}, extending \cite{bartsch})
that a necessary condition for the existence of a 
$\Zz /p^{k+1}$-map $S(nL)\to S(V)$ is: $n\leq \delta (V)$.
\Thm{\label{pos}
There is a constant $c\geq p^k-1$, depending on $p$ and $k$,
such that if $\sum_{l=0}^k m_lp^l >c$ there exists
a $\Zz /p^{k+1}$-equivariant map 
$f_0: S(((\sum m_lp^l)-c)L) \to S(V)$.
For $n> \sum m_lp^l -c$, the map
$$\textstyle
f : S(nL)=S((n-(\sum m_lp^l) +c)L)*S(((\sum m_lp^l)-c)L) \to V
$$ 
defined by $f([u,t,v]) = tf_0(v)$ has zero-set
the sphere $S((n-\sum m_lp^l+c)L)$,
of dimension $2(n-\sum m_lp^l+c)-1$.
}

For results in this direction, requiring substantial input
from stable homotopy theory \cite{SS, meyer},
see \cite[Section 4]{bartsch2} and \cite[Theorem 5.5]{meyer}. 

Already in \cite{nakaoka} the vector space $V$ in Theorem \ref{main}
was replaced, in certain cases, by a smooth manifold.
Suppose that $M$ is a connected smooth $\Zz /p^{k+1}$-manifold (without boundary and
admitting an embedding as a $\Zz /p^{k+1}$-submanifold in some Euclidean 
$\Zz /p^{k+1}$-module) with non-empty fixed submanifold 
$M^{\Zz /p^{k+1}}$.
\Thm{\label{mainX}
Let $X$ be a compact connected free $\Zz /p^{k+1}$-ENR such that
$H_i(X;\,\Ff_p)$ is zero in dimensions $0<i<2(n-1)$.

Consider a $\Zz /p^{k+1}$-equivariant map
$f : X \to M$
which is equivariantly homotopic to the constant map 
at some point $z_0$ in the fixed submanifold $M^{\Zz /p^{k+1}}$.

Let $N\subseteq M$ denote the component of the fixed submanifold
containing $z_0$, and suppose that the normal space
$\tau_{z_0}M/\tau_{z_0}N$ at $z_0$ is isomorphic
as a real $\Zz /p^{k+1}$-module to the complex $\Zz /p^{k+1}$-module
$V$ satisfying the conditions that
$m_k\not=0$ and $n>\delta (V)$.
If $p=2$, we make the additional assumption that $m_0=0$.

Then there exists
an integer $d\geq 2(n-1-\delta (V))$ such that 
the mod $p$ cohomology group
$H^d(f^{-1}(N);\, \Ff_p)$ of the compact set
$$
f^{-1}(N) = \{ x \in X \st f(x)\in N\}
$$
is non-trivial.
}
If $p=2$, but $m_0\not=0$, the conclusion remains true
if the $2m_0$-dimensional real subbundle of $\tau M/\tau N$ over $N$ on
which the generator of $\Zz /2^{k+1}$ acts as $-1$ admits
a complex structure.
\Cor{\label{quotient}
Suppose, further, under the hypotheses of Theorem \ref{mainX}
that a finite group $\Gamma$ acts on $X$,
commuting with the action of $\Zz /p^{k+1}$, and that
$f$ is the composition of the projection $X\to X/\Gamma =\bar X$ with
a $\Zz /p^{k+1}$-map $\bar f : \bar X \to M$. Then the subspace
$$
\{ \bar x\in X/\Gamma \st \bar f(\bar x) \in N\}
$$
has covering dimension greater than or equal to
$2(n-1-\delta (V))$. 
}

The Borsuk--Ulam theorem for connective $K$-theory
is recalled in Section 2. Theorem \ref{main}
and Corollary \ref{quotient} are established in Section 3,
and the constructions required to prove Theorem \ref{pos}
are given in Section 6.

This note should be read as
a supplement to \cite{kbu}, explaining how results
presented in that paper for the prime $2$ can be adapted to the 
case of an arbitrary prime.
\Sect{Connective $K$-theory Borsuk--Ulam theorems}
The spaces that one meets in classical Algebraic Topology are 
usually CW-complexes. 
Here we shall work in the category
of compact Hausdorff spaces (or, in fact, compact subspaces of
Euclidean spaces).  To extend standard results to these spaces
it is natural to use representable cohomology theories,
as described, for example, in \cite[Section 8]{borsuk}.
Connective complex $K$-theory, $k^*$, as well as $\Ff_p$-cohomology,
is thus to be understood as a representable theory, represented
by spaces which are CW-complexes (thought of as unions of finite
complexes). 

An $m$-dimensional complex vector bundle $\xi$ over a compact ENR 
(Euclidean Neighbourhood Retract) $Y$
has a connective $K$-theory Euler class $e_k(\xi )\in k^{2m}(Y)$
(defined as the restriction to the zero section of the Bott class 
of $\xi$ in the $K$-theory
$k^{2m}(D(\xi ),S(\xi ))$ of the disc bundle $D(\xi)$ modulo the
sphere bundle $S(\xi )$).

There is a discussion of the following key result in \cite{kbu};
see, also, \cite[Proposition 2.7]{borsuk}. A proof of a more
general statement, Theorem \ref{bloc}, is given later.
\Thm{\label{loc}
{\rm (A connective $K$-theory Borsuk--Ulam theorem).}
Let $\xi$ be an $m$-dimensional complex vector bundle over a compact ENR $Y$, and
let $s$ be a section of $\xi$ with zero-set $\Zero (s)$.
Suppose that $A\subseteq Y$ is a closed sub-ENR
and that $a\in k^i(Y,A)$ is a class such that 
$a\cdot e_k(\xi )\not=0\in k^{i+2m}(Y,A)$. 
Then $a$ restricts to a non-zero class
in $k^i(\Zero (s), \Zero (s)\cap A)$. 
\qed
}
We can also look, for $l \geq 1$,
at (mod $p^l$) connective $K$-theory defined
by
$$
(k\Zz /p^l)^i(Y,A) = k^{i+2}((Y,A) \times (C_l,*)),
$$
where $C_l$ is the cofibre of the $p^l$-th power map
$$
S^1 \Rarr{z\mapsto z^{p^l}}{} S^1 \to C_l
$$
on $S^1=S(\Cc )$ with basepoint $*=1$.
There is a Bockstein exact sequence
$$
\cdots\Rarr{p^l}{}
k^i(Y,A)=k^{i+1}((Y,A)\times (S^1,*))
\to
(k\Zz /p^l)^i(Y,A)
$$
$$
=k^{i+2}((Y,A)\times (C_l,*))
\to
k^{i+1}(Y,A) =k^{i+2}((Y,A)\times (S^1,*))
\to\cdots \, .
$$

\smallskip

The mod $p^l$ Borsuk--Ulam theorem follows at once from 
Theorem \ref{loc}.
\Prop{\label{ploc}
Let $\xi$ be an $m$-dimensional complex vector bundle over a compact ENR $Y$.
Let $s$ be a section of $\xi$ with zero-set $\Zero (s)$.
Suppose that $A\subseteq Y$ is a closed sub-ENR
and that $a\in (k\Zz /p^l)^i(Y,A)$ is a cohomology class 
such that $a\cdot e_k(\xi )\not=0\in (k\Zz /p^l)^{i+2m}(Y,A)$.
Then $a$ restricts
to a non-zero class in $(k\Zz /p^l)^i(\Zero (s), \Zero (s)\cap A)$. 
\qed
}
Consider more generally, for some compact Lie group $G$,
a smooth $G$-manifold $M$ (finite-dimensional,
without boundary, and admitting an embedding
as a submanifold of some Euclidean $G$-module) and a 
$G$-submanifold $N\subseteq M$ without boundary, embedded
as a closed subspace of $M$. 
Let $P\to Y$ be a principal $G$-bundle over a compact ENR $Y$.
Now form the fibre bundle $E=P\times_G M\to Y$ and the subbundle
$F=P\times_GN\subseteq E$ over $Y$. We write $\nu$ for the fibrewise
normal bundle over $F$ of the inclusion $F\subseteq E$:
the restriction of $\nu$ to the fibre $F_y$ over a point $y\in Y$
is the normal bundle of the submanifold $F_y\subseteq E_y$. 
\Thm{\label{bloc}
In the notation introduced above,
suppose that $z$ is a section of the subbundle $F\to Y$ and that
$s$ is a section of the fibre bundle  $E\to Y$ 
homotopic to the section $z$.
Write
$$
\Null (s) =\{ y\in Y \st s(y) \in F\}.
$$
Assume that the fibrewise normal bundle $\nu$ admits the structure
of a complex vector bundle of dimension $m$, and
write $\xi$ for the complex vector bundle $z^*\nu$ over $Y$.
Suppose that we are given a class $a\in k^i(Y)$. 
If $a\cdot e_k(\xi )\not=0\in k^{i+2m}(Y)$, 
then the restriction of $a$ to the compact subspace
$\Null (s)$ of $Y$ is a non-zero class in $k^i(\Null (s))$.
}
The construction that will be used in the proof is taken from
coincidence theory, \cite[Definition 7.2]{coin};
but the ideas are implicit in \cite{nakaoka}.
\begin{proof}
We shall prove that if $a$ restricts to zero in
$k^i(\Null (s))$, then $a\cdot e_k(\xi )=0$.

So, assuming that $a$ restricts to zero on $\Null (s)$,
we can choose an open neighbourhood $U$ of $\Null (s)$ in $Y$
such that $a$ restricts to zero in $k^i(\cl{U})$.
Thus $a$ lifts to some class $a'\in k^i(Y,\cl{U})$.
We shall show that $e_k(\xi )$ lifts to a class
$e'\in k^{2m}(Y,Y-U)$. Then, because $Y$ is the union of the closed 
subspaces $\cl{U}$ and $Y-U$, $a'\cdot e'=0$, and hence
$a\cdot e_k(\xi )=0$.

Fix a homotopy $H : [0,1]\times Y\to E$ between
$z$ and $s$. Since $Y$ is compact,
there is a compact $G$-subspace $K$ of $M$ such that
$H([0,1]\times Y)\subseteq P\times_G K$.
Next, choose open $G$-subspaces $W, \, W'\subseteq M$ such
that $\cl{W}\subseteq W'$, $\cl{W'}$ is compact and $K\subseteq W$.
(If $M$ is compact, we can take $W=W'=M$.)
Set $F_0=P\times_G (F\cap W)$ and $\cl{F_0}=P\times_G (F\cap \cl{W})$,
$F_0'=P\times_G (F\cap W')$ and $\cl{F_0'}=P\times_G (F\cap \cl{W'})$.
Choose a $G$-invariant Riemannian metric on $M$. Then, for
$\epsilon >0$ sufficiently small, the exponential map gives
a fibrewise tubular neighbourhood embedding the closed disc bundle
$D_\epsilon (\nu \, |\, \cl{F_0'})$ of radius $\epsilon$
in $E$ and the open disc bundle $B_\epsilon (\nu \, |\, F_0')$ 
as an open neighbourhood of $F'_0$ with the property that
$D_\epsilon (\nu\, |\, \cl{F_0'}-F_0')$ is disjoint from the
compact subspace $P\times_G K$.
Since $s(y)\notin F$ if $y\in Y-U$, by reducing $\epsilon$ if
necessary, we can arrange that 
$s(y)\notin D_\epsilon (\nu\, |\, F_0)$ for $y\in Y-U$.
We can now define a continuous map 
$$
\sigma: Y\to D_\epsilon (\nu\, |\, \cl{F_0})/
S_\epsilon (\nu\, |\, \cl{F_0}),
$$
sending $y\in Y$ to $[s(y)]$ if $s(y)\in D_\epsilon (\nu \, |\, Y_0)$
and otherwise to the basepoint $*$.
Since $\sigma$ maps $Y-U$ to $*$, we have a homomorphism
$$
\sigma^* : k^{2m}(D_\epsilon (\nu\, |\, \cl{F_0}),
S_\epsilon (\nu\, | \,\cl{F_0}))\to k^{2m}(Y,Y-U).
$$
Now define $e'\in k^{2m}(Y,Y-U)$ to be the image
of the Bott class of the complex bundle $\nu$.
Applying the same construction to the section $z$,
we get a map
$$
\zeta: Y\to D_\epsilon (\nu\, |\, \cl{F_0})/
S_\epsilon (\nu\, |\, \cl{F_0}),
$$
which is homotopic to $\sigma$. But $\zeta$ is just the inclusion
of $Y$ into the zero-section 
$\cl{F_0}\into D_\epsilon (\nu\, |\, \cl{F_0})$, and so 
$\zeta^*$ maps the Bott class, by definition, to the Euler
class $e_k(\xi )$. Thus $e'$ lifts $e_k(\xi )$, as required.
\end{proof}

\Sect{A lower bound on the dimension of the zero-set}
Let $X$ be a compact connected free $\Zz /p^{k+1}$-ENR with
$H_i(X;\,\Ff_p)=0$ for $0<i<2(n-1)$.

Consider first a $\Zz /p^{k+1}$-equivariant map $f: X\to V$.
Form the complex vector bundle $\xi = (X\times V)/(\Zz /p^{k+1})$
over the orbit space $Y=X/(\Zz /p^{k+1})$ and the complex line
bundle $\lambda = (X\times L)/(\Zz /p^{k+1})$
(associated with the basic $1$-dimensional representation
$L$ of $\Zz /p^{k+1}$).
The map $f$ determines a section $s$ of $\xi$:
$$
s([x])=[x,f(x)].
$$
Let us write $Z=\Zero (s)\subseteq Y$ and $\tilde Z=\{ x\in X\st
f(x)=0\}$. Then $Z=\tilde Z/(\Zz /p^{k+1})$.

\Lem{\label{appl}
Suppose that $e_k(\lambda)^j\cdot e_k(\xi)\in k^{2j+2m}(Y)$
is non-zero.
Then the restriction of $e_k(\lambda )^j$ to $k^{2j}(Z)$ is non-zero 
and, moreover, if $l\geq 1$ is sufficiently large
the image of $e_k(\lambda)^j$ in $(k\Zz /p^l)^{2j}(Z)$ is non-zero.
}
In the  statement $e_k(\lambda )^0$ is to be read as $1$.
\begin{proof}
The restriction of $e_k(\lambda )^j$ to $k^{2j}(Z)$ is non-zero
by Theorem \ref{loc}.

The $k$-theory Euler class, $e_k(\lambda )$, of $\lambda$
is a $p$-primary torsion class, because its lift to the 
$p^{k+1}$-fold cover $X$ of $Y$ is zero.
(Consider the projection $\pi : X\to Y$ and the
transfer $\pi_!: k^*(X)\to k^*(Y)$. We see that
$\pi_!(1)\cdot e_k(\lambda )=\pi_!(\pi^*e_k(\lambda ))=
\pi_!(e_k(\pi^*\lambda ))=0$ and 
$\pi_!(1)-p^{k+1}\in k^0(Y)$ is nilpotent.).
Since the group $k^{2j+m}(Y)$ is finitely generated, for $l\geq 1$
sufficiently large the (non-zero) $p$-primary torsion class
$e_k(\lambda )^je_k(\xi)$ is not divisible by $p^l$,
and hence the image of 
$e_k(\lambda )^je_k(\xi )$ in the mod $p^l$ group
$(k\Zz /p^l)^{2j+m}(Y)$
is non-zero. We conclude from the mod $p^l$ Borsuk--Ulam theorem,
Proposition \ref{ploc}, that the image of $e_k(\lambda )^j$ in 
$(k\Zz /p^l)^{2j}(Z)$ is non-zero.
\end{proof}

We perform the calculation for a specific space $X$. 
Choose a triangulation of $S(nL)/(\Zz /p^{k+1})$ 
as a $(2n-1)$-dimensional
finite complex, and let $Y_0$ be its $2(n-1)$-skeleton.
Take $X_0$ to be the inverse image of $Y_0$ in
$S(nL)$ and write $\xi_0$ and $\lambda_0$ for the corresponding
complex vector bundles over $Y_0$.
\Lem{\label{special}
Suppose that $m_k\not=0$ and $n>\delta (V)$.
Then, if $0\leq j\leq n-1-\delta (V)$,
$e_k(\lambda_0)^j\cdot e_k(\xi_0)\in k^{2j+2m}(Y_0)$
is non-zero.
}
\begin{proof}
It is enough to verify that the image of $e_k(\lambda_0)^je_k(\xi_0 )$
in the simpler periodic $K$-theory is non-zero.

As a real $\Zz /p^{k+1}$-module, $V$ is isomorphic to
$$
\bigoplus_{i=1}^m L^{\otimes t_i},\quad
\text{where $t_i\geq 1$ and $m_l =\#\{ i\st \nu_p(t_i)=l\}$,}
$$
and we may assume that $V$ is precisely this module.
So $\xi_0=\bigoplus_{i=1}^m \lambda_0^{\otimes t_i}$.

The quotient $S(nL)/\Tt$ by the action of the group $\Tt$ of complex numbers of modulus $1$ is the $(n-1)$-dimensional complex
projective space $\Cc P(nL) =\Cc P(\Cc^n)$.
We write $H$ for the (dual) Hopf line bundle over
$\Cc P(nL)$ with fibre $\Cc v$ over $[v]$, ($v\in S(nL)$).

We can now identify $S(nL)/(\Zz /p^{k+1}\Zz)$ with the sphere bundle
$S(H^{\otimes p^{k+1}})$ of the $p^{k+1}$-th power of
$H$ over 
$\Cc P(nL)$ by mapping $v\in S(nL)$ to $v^{p^{k+1}}$ in the fibre
of $H^{\otimes p^{k+1}}$ over $[v]\in \Cc P(nL)$.

Let $z=[H]\in K^0(\Cc P(\Cc^n))$. 
Recall that $K^0(\Cc P(\Cc^n))=\Zz [z]/(1-z)^n$
and $K^1(\Cc P(\Cc^n))=0$.
The $K$-theory Euler class of $H$ and, more generally,
of $H^{\otimes t}$ can be written,
using the periodicity, as $e_K(H)=1-z\in K^0(\Cc P(\Cc^n))$
and $e_K(H^{\otimes t})=1-z^t$.

The group $K^0(S(H^{\otimes p^{k+1}}))$ is calculated from 
the Gysin sequence of the sphere bundle
$$
K^0(\Cc P(\Cc^n)) \Rarr{1-z^{p^{k+1}}}{}
K^0(\Cc P(\Cc^n)) \to K^0(S(H^{\otimes p^{k+1}})) \to 0
$$
as $\Zz [z]/(1-z^{p^{k+1}},\, (1-z)^n)$.
(Another approach would use equivariant $K$-theory.
This is the quotient of the representation ring
$K^0_{\Zz /p^{k+1}}(*)=R(\Zz /p^{k+1})$
by the ideal generated by the $\Zz /p^{k+1}$-Euler class of $nL$.)

We shall prove in Section \ref{algebra} that the class of
the polynomial 
$$
(\prod_{i=1}^m(1-z^{t_i}))(1-z)^j \in 
\Zz [z]/(1-z^{p^{k+1}},\, (1-z)^n)
$$ 
is non-zero, that is, $e_K(H)^j\cdot e_K(\bigoplus_{i=1}^m H^{\otimes t_i}) \not=0$ in $K^0(S(H^{\otimes p^{k+1}}))$.

Now the restriction
$$
K^0(S(nL)/(\Zz /p^{k+1})) \to K^0(Y_0)
$$
is injective, because $K^0(D(\Rr^{2n-1}),S(\Rr^{2n-1}))=0$.
So $e_K(\lambda_0 )^j\cdot e_k(\xi_0 )\not=0$.
\end{proof}
The next technical lemma will be proved in Section \ref{technical}.
\Lem{\label{reduction}
If $e_k(\lambda_0)^j\cdot e_k(\xi_0)\in k^{2j+2m}(Y_0)$ is non-zero,
then $e_k(\lambda)^j\cdot e_k(\xi )\in k^{2j+2m}(Y)$ is non-zero.
}
In conjunction with Lemmas \ref{appl} and \ref{special} this 
establishes the main result.
\Prop{\label{appln}
Suppose that $m_k\not=0$ and $n>\delta (V)$.
Then, if $0\leq j\leq n-1-\delta (V)$, 
the group $k^{2j}(Z)$ is non-zero 
and, if $l\geq 1$ is sufficiently large,
$(k\Zz /p^l)^{2j}(Z)$ is non-zero.
\qed
}
This implies bounds on the covering dimension and cohomological
dimension of he zero-set $Z$.
\Cor{\label{cohdim}
Let $d_0= 2(n-1-\delta (V))$.
Then 
\par\noindent {\rm (i)}\ \ \ 
the covering dimension of $Z$ is at least $d_0$;
\par\noindent {\rm (ii)}\ \ 
there is an integer $d_1\geq d_0$ such that
$H^{d_1}(Z;\, \Zz)\not=0$;
\par\noindent {\rm (iii)}\ 
there is an integer $d_2\geq d_0$ such that
$H^{d_2}(Z;\, \Ff_p)\not=0$.
}
\begin{proof}
According to Proposition \ref{appln}, the groups
$k^{d_0}(Z)$ and, for some $l\geq 1$, $(k\Zz /p^l)^{d_0}(Z)$
are non-zero.

Let $d_1\geq d_0$ be the greatest integer such that $k^{d_1}(Z)$
is non-zero; by Lemma \ref{bound}, we have $d_1\leq 2n-1$.
And let $d_2\, (\,\leq d_1)$ be the greatest integer such that
$(k\Zz /p^l)^{d_2}(Z)$ is non-zero.

The assertions (i), (ii) and (iii) then
follow from Lemmas \ref{cdim}, \ref{coho} and \ref{coho2} in
the Appendix.
\end{proof}
To establish the corresponding result for $\tilde Z$ we
consider the tower
$$
\tilde Z \to \tilde Z/(p^k\Zz /p^{k+1})
\to\tilde Z/(p^{k-1}\Zz /p^{k+1})\to \cdots\to
\tilde Z /(\Zz /p^{k+1}) =Z
$$
of principal $\Zz /p$-bundles.
\Cor{\label{cohdim2}
There exist integers $d_{\tilde 1}$ and $d_{\tilde 2}$ greater
than or equal to $d_0$ such that the groups
$k^{d_{\tilde 1}}(\tilde Z)$,
$(k\Zz /p^l)^{d_{\tilde 2}}(\tilde Z)$,
$H^{d_{\tilde 1}}(\tilde Z;\,\Zz )$
and 
$H^{d_{\tilde 2}}(\tilde Z;\, \Ff_p)$
are all non-zero.
The covering dimension of $\tilde Z$ is at least $d_0$.
}
\begin{proof}
Successive application of Lemma \ref{cover}(ii) to the 
$k+1$ principal $\Zz /p$-bundles
shows that there is some $d\geq d_0$ such that
$k^d(\tilde Z)$ is non-zero. Let $d_{\tilde 1}$,
bounded like $d_1$ by $2n-1$, be
the greatest such integer $d$ and let $d_{\tilde 2}$ be the 
greatest integer such that $(k\Zz /p^l)^{d_{\tilde 2}}(\tilde Z)$
is non-zero.
Now apply Lemmas \ref{coho} and \ref{coho2} again,
and then Lemma \ref{cdim}, to conclude that the covering
dimension is at least $d_{\tilde 2}$.
\end{proof}

This completes the proof of Theorem \ref{main}
with $d=d_{\tilde 2}$.

\begin{proof}[Proof of Theorem \ref{mainX}]
We apply Theorem \ref{bloc} with $G=\Zz /p^{k+1}$ acting, as
given, on
$M$ and (trivially) on $N$ and with the quotient $X\to X/G$
as the principal bundle $P\to Y$.
The map $f$ determines the section $s$:
$s([x])=[x,f(x)]$, $x\in X$,
and the point $z_0$ determines
the section $z([x])=[x,z_0]$.
Every irreducible real representation of $\Zz /p^{k+1}$ is
complex, with the exception of the $1$-dimensional representation
if $p=2$.
The normal bundle $\nu$, therefore, admits a complex structure
and the vector bundle $\xi=z^*\nu$ is isomorphic to 
$X\times_{\Zz /p^{k+1}} V$.

Take $a=e_k(\lambda )^j$, where $j=n-1-\delta (V)$,
and apply Lemmas \ref{appl}, \ref{special}
and \ref{reduction}.
\end{proof}
\Cor{\label{mn}
Let $X$ be a compact connected free $\Zz /p^{k+1}$-ENR such that
the $\Ff_p$-cohomology group
$H^i(X;\,\Ff_p)$ is zero in dimensions $0<i<2(n-1)$.
Let $N$ be a connected smooth manifold of dimension $r$. 

Consider a map $h : X \to N$
which is homotopic to a constant map.

Assume that $p$ is odd.
Then there exists
an integer 
$$
d\geq 2(n+(p^k-1)-1)-r(k+1)(p-1)p^k\, ,
$$
such that the mod $p$ cohomology group
$H^d(A(h );\, \Ff_p)$ of the compact set
$$
A(h ) = \{ x \in X \st h(\gamma x)=h(x)
\text{\ for all $\gamma\in\Zz /p^{k+1}$}\}
$$
is non-trivial.
}
If $p=2$, the conclusion holds under the additional assumption
that the tangent bundle $\tau N$ admits a complex structure.
\begin{proof}
Put $M=\map (\Zz /p^{k+1},N)$ with the permutation action of
$\Zz /p^{k+1}$ and include $N$ in $M$
as the $\Zz /p^{k+1}$-fixed subspace of constant maps.
Let $f: X\to M$ be the $\Zz /p^{k+1}$-map
given by $f(x)(\gamma )=h(\gamma x)$.
Thus $A(h)=\{ x\in X \st f(x)\in N\}$.
Choose a point $z_0\in N$. Since $h$ is homotopic
to the constant map $z_0$, the map $f$ is
equivariantly homotopic to $z_0$.

We now apply Theorem \ref{mainX}.
If $p$ is odd, the normal
space at $z_0\in N\subseteq M$ is isomorphic to
$V$ with $m_l =r(p-1)p^{k-l}/2$, so that
$\sum_{l=0}^k p^lm_l  = (r(k+1)(p-1)/2)p^k$
and $m_k\not=0$.
And the normal bundle of the inclusion $N\subseteq M$
admits an equivariant complex structure.

If $p=2$, the additional assumption ensures that $r$
is even and that the normal bundle admits a complex structure.
\end{proof}
\Rem{In \cite[Theorem 2]{munkholm}, $p$ is odd and
$\Sigma$ is a homotopy sphere of dimension $2n-1$ with
a free action of $\Zz /p^{k+1}$. Assuming that $\Sigma$ is a smooth
$\Zz /p^{k+1}$-manifold, we can take $X$ to be the complement of
a small equivariant open tubular neighbourhood of the orbit of
some point of $\Sigma$. The restriction to $X$ of any map
defined on $\Sigma$ is null-homotopic. So we may apply
Corollary \ref{mn}.
} 
\begin{proof}[Proof of Corollary \ref{quotient}]
Composing the projection $a: X\to \bar X$ with $\bar f$ we get a 
$\Zz /p^{k+1}$-map 
$$
f= \bar f\comp a : X \to V\, . 
$$
Then $\tilde Z =\{ x\in X \st f(x)=0\} = a^{-1}(\bar Z)$
and $\tilde Z \to \bar Z$ is the quotient by
the action of $\Gamma$.
The lower bound on the covering dimension of $\bar Z =\tilde Z/\Gamma$ 
now
follows from Proposition \ref{orbit} and Corollary \ref{cohdim2}.
\end{proof}
\Cor{\label{general}
Let $U$ be a complex $\Zz /p^{k+1}$-module of dimension $n$
and consider a $\Zz /p^{k+1}$-map
$$
\bar f : S(U) \to V.
$$

If $m_k\not=0$ and $n>\delta (V)$,
then the covering dimension of the zero-set
$\bar Z=\{ v\in S(U) \st \bar f (v)=0\}$
is at least $2(n-\delta (V)-1)$.
}
\begin{proof}
Suppose that, as a real module,
$U=\bigoplus_{i=1}^n L^{\otimes r_i}$, where $1\leq r_i\leq k+1$.
We have a $\Zz /p^{k+1}$-map
$$
a: S(nL)=S(L\otimes_\Cc \Cc^n) \to S(U)\ :\ (z_1,\ldots ,z_n)\mapsto
(z_1^{r_1},\ldots ,z_n^{r_n})/(\sum_i |z_i|^{2r_i})^{1/2}\, .
$$
Take $X=S(nL)$ and $\bar X=S(U)$, so that $\bar X$ is the
quotient of $X$ by the action of
a finite group $\Gamma$ (the product of the groups of $r_i$th roots
of unity) of order $r_1\cdots r_n$.
The result follows from Corollary \ref{quotient}.
\end{proof}
\Rem{This result strengthens \cite[Theorem 4.2]{BMS}, 
which shows that, if, for some $k'\leq k$,
$V^{p^{k-k'}\Zz /p^{k+1}}=0$ (so that $m_l=0$ for $l>k'$) and
$U^{p^{k-k'}\Zz /p^{k+1}}=0$, then 
$\dim \bar Z \geq 2r$, where $r$ is the smallest integer
such that $r\geq (n-1-p^{k'}m)/p^{k'}$. 

It also, when $p=2$
and $m_k=m$, strengthens \cite[Theorem 1.1]{kbu},
which gives the bound $2(n-2^km-1)$,
rather than $2(n-2^km+(2^k-1)-1)$.
}
\Sect{\label{algebra}
An algebraic lemma}
A proof of the basic algebraic lemma is implicit in the work of
Ohashi \cite[Theorem 2.3]{ko}, generalizing special cases in
\cite[Section 2]{munkholm},
\cite[Proposition 3.1]{bartsch}. The proofs in \cite{munkholm, ko}
use results about cyclotomic polynomials.
For the sake of completeness we include an elementary direct proof.
(More detailed descriptions of the $K$-theory of the 
lens space $S(nL)/(\Zz /p^{k+1})$ can be found in the lierature, 
for example, in \cite{kobayashi}.)

Recall that we are taking $V=\bigoplus_{i=1}^m L^{\otimes t_i}$,
where $t_i\geq 1$, and assuming that $m_k\not=0$.
\Lem{{\rm (\cite[Theorem 2.3]{ko}).}
The class of $(\prod_{i=1}^m(1-z^{t_i}))(1-z)^j$ 
is non-zero in 
$\Zz [z]/(1-z^{p^{k+1}},\, (1-z)^n)$
if $n\geq j+1+\delta (V)$.
}
\begin{proof}
Let $\phi (z)=1+z^{p^k}+\ldots +z^{(p-1)p^k}\in\Zz [z]$,
so that $1-z^{p^{k+1}}=(1-z^{p^k})\phi (z)$.
 
Calculating modulo $p$, we have, for $0\leq l\leq k$, 
$$
(1-z^{p^l})^{(p-1)p^{k-l}} \equiv(1-z^{p^k})^{p-1}
\equiv (1-z^{p^k})^p(1-z^{p^k})^{-1}
\equiv (1-z^{p^{k+1}})(1-z^{p^k})^{-1}
$$
and so can write 
$$
(1-z^{p^l})^{(p-1)p^{k-l}} = - p(1+(1-z)a_l(z))+ \phi (z),
\leqno{{\rm (a)}}
$$
where $a_l(z)\in\Zz [z]$, since $\phi (1)=p$.

We now localize at the prime $p$ and work in the Unique Factorization
Domain $\Zz_{(p)}[z]$.

Suppose that $t\geq 1$ is not divisible by $p$. Then
$1-z^{tp^l}=(1-z^{p^l})(t+(1-z)h(z))$, where $h(z)\in \Zz [z]$,
and $t+(1-z)h(z))$ is a unit in $\Zz_{(p)}[z]/((1-z)^n)$.
It suffices, therefore, to show that
$$
(\prod_{l=0}^k (1-z^{p^l})^{m_l})(1-z)^j\in
\Zz_{(p)}[z]/(1-z^{p^{k+1}},\, (1-z)^n)
$$
is non-zero. 

Write $m'_k=m_k-1$ and $m'_l=m_l$ for $l<k$
and set
$$
E(z)=(\prod_l(1-z^{p^l})^{m'_l})(1-z)^j.
$$ 
We have to show that $(1-z^{p^k}) E(z)$ 
is non-zero in the ring
$$
\Zz_{(p)} [z]/((1-z^{p^k})\phi(z),\, (1-z)^{j+\delta (V)+1}),
$$
or, equivalently, that
$E(z)$ is non-zero in the quotient ring
$$
\Zz_{(p)} [z]/(\phi(z),\, (1-z)^{j+\sum m'_lp^l +1})=
\Zz_{(p)} [z]/(\phi(z),\, (1-z)^{j+\delta (V)}),
$$
because $1-z$ divides $1-z^{p^k}$, but $(1-z)^2$ does not.

It follows from (a) that, for each $l\leq k$,
because $1+(1-z)a_l(z)$ is a unit in $\Zz [z]/(1-z)^n$,
$p$ lies in the principal ideal in
$\Zz_{(p)}[z]/(\phi (z), (1-z)^n)$ generated by $1-z^{p^l}$.
We also have $1-z^{p^l}\equiv (1-z)^{p^l}$ (mod $p$).  Hence,
for $0\leq l\leq k$,
$(1-z)^{p^l}$ lies in the ideal generated by $1-z^{p^l}$.

So $(1-z)^{j+\sum m'_lp^l}$ lies in the ideal generated by $E(z)$.
Since this power is non-zero in the quotient ring, so is $E(z)$.
\end{proof}
\Rem{
Using (a) with $l=0$, we see that $p$ lies in the principal ideal in
$\Zz_{(p)}[z]/(\phi (z), (1-z)^n)$ generated by $(1-z)^{(p-1)p^k}$
and hence that,
for $0\leq l\leq k$,
$1-z^{p^l}$ (which is congruent (mod $p$) to $(1-z)^{p^l}$)
lies in the ideal generated by $(1-z)^{p^l}$.
So $E(z)$ lies in the ideal generated by
$(1-z)^{\sum m'_lp^l +j}$. The bound on $n$ is, thus, optimal.
}
\Sect{\label{technical}
Reduction to the case of a sphere}
\Lem{\label{Htok}
Suppose that $X$ is compact, connected ENR such that
$H_i(X;\,\Ff_p)=0$ for $0<i<2(n-1)$. 
Then $\tilde k_i(S^0*X)_{(p)}=0$ for $i\leq 2(n-1)$. 
}
\begin{proof}
We have a cofibre sequence
$$
X_+ \to S^0 \to S^0*X,
$$
and so a long exact sequence
$$
\cdots\to
H_i(X;\,\Ff_p) \to H_i(*;\,\Ff_p) \Rarr{0}{} \tilde H_i(S^0*X;\,\Ff_p)
\to H_{i-1}(X;\,\Ff_p)\to\cdots
$$
Thus $\tilde H_i(S^0*X;\,\Ff_p)=0$ for $i\leq 2(n-1)$.

Since $\tilde H_i(S^0*X;\,\Zz )$ is a finitely generated abelian group,
we deduce from the Bockstein exact sequence
$$
\to \tilde H_{i+1}(S^0*X;\,\Ff_p)\to
\tilde H_i(S^0*X;\, \Zz )\Rarr{p}{} \tilde H_i(S^0*X;\, \Zz )\to
\tilde H_i(S^0*X;\, \Ff_p)\to
$$
that $\tilde H_i(S^0*X;\,\Zz )_{(p)}=0$ for $i\leq 2(n-1)$.

From the exact sequence 
$$
 \to \tilde H_{i+1}(S^0*X;\Zz)_{(p)}\to
\tilde k_{i-2}(S^0*X)_{(p)}\Rarr{v\cdot}{} \tilde k_i(S^0*X)_{(p)}
\to \tilde H_i(S^0*X;\, \Zz)_{(p)}\to
$$
relating connective $K$-theory and homology by multiplication by the Bott generator $v\in k_2(*)=\Zz v$, we see that
$\tilde k_i(S^0*X)_{(p)}=0$ for $i\leq 2(n-1)$. 
\end{proof}
\Lem{\label{split}
Under the hypotheses of Lemma \ref{Htok} the projection
$$
\pi_0 : (X_0\times X)/(\Zz /p^{k+1})\to X_0/(\Zz /p^{k+1})=Y_0
$$
to the first factor induces a split injection
$$
k^* (Y_0)_{(p)} \to k^*((X_0\times X)/(\Zz /p^{k+1}))_{(p)}\, .
$$
}
\begin{proof}
We shall use fibrewise connective $K$-theory, as in 
\cite[II, Proposition 15.42]{fht},
applied to the fibrewise cofibre sequence over $Y_0$:
$$
(X_0\times X_+)/(\Zz /p^{k+1})\to
(X_0\times S^0)/(\Zz /p^{k+1})\to
(X_0\times (S^0*X))/(\Zz /p^{k+1})
$$
in which the space $X_+$ obtained by adding a disjoint basepoint 
to $X$ projects to $\{ *\}_+=S^0$.
From Lemma \ref{Htok}, by induction over the cells of
the $2(n-1)$-dimensional complex $Y_0$, we see that
$$
k^0_{Y_0}\{ Y_0\times S^0;\, 
(X_0\times (S^0*X))/(\Zz /p^{k+1})\}_{(p)}
=0\, .
$$
So
$$
k^0_{Y_0}\{ Y_0\times S^0;\, (X_0\times X_+)/(\Zz /p^{k+1})\}_{(p)}
\to k^0_{Y_0}\{ Y_0\times S^0;\, Y_0\times S^0\}_{(p)}
$$
is surjective,
and there is an element in
$$
k^0_{Y_0}\{ Y_0\times S^0;\, (X_0\times X_+)/(\Zz /p^{k+1})\}_{(p)}
$$
lifting the identity 
$1\in k^0(Y_0)_{(p)}
=k^0_{Y_0}\{ Y_0\times S^0;\, Y_0\times S^0\}_{(p)}$.
Collapsing the fibrewise basepoints to a point, we get an element
$\iota\in k^0\{ (Y_0)_+;\, ((X_0\times X)/(\Zz /p^{k+1}))_+\}_{(p)}$
such that $(\pi_0)_* \comp\iota$ is the identity in
$k^0\{ (Y_0)_+;\, (Y_0)_+\}_{(p)}$.

This element $\iota$ determines the required splitting.
\end{proof}
\Rem{Alternatively, we could use stable homotopy theyr,
written as $\omega_*$. By the Whitehead theorem, 
$\tilde\omega_i(S^0* X)_{(p)}=0$ for $i\leq 2(n-1)$.
Running the argument in Lemma \ref{split} for stable
homotopy, we get an element in
$$
\omega^0_{Y_0}\{ Y_0\times S^0;\, 
(X_0\times X_+)/(\Zz /p^{k+1})\}_{(p)},
$$
which can be thought of as a $p$-local stable cross-section
of $\pi_0$.
}
\begin{proof}[Proof of Lemma \ref{reduction}]
Consider the projection
$$
\pi : (X_0\times X)/(\Zz /p^{k+1})\to X/(\Zz /p^{k+1})=Y
$$
to the second factor.
Now
$\pi_0^*(e_k(\lambda_0)^j\cdot e_k(\xi_0))$ is equal to
$\pi^*(e_k(\lambda )^j\cdot e_k(\xi ))$.
It, therefore, follows from the injectivity of $\pi_0^*$
established in Lemma \ref{split}
that, if the class $e_k(\lambda_0)^j\cdot e_k(\xi_0)$,
which we know to be $p$-primary torsion, is non-zero,
then so, too, is $e_k(\lambda )^j\cdot e_k(\xi )$.
\end{proof}
\Sect{Constructions}
For $t\geq 0$, we identify the $t$th complex tensor power 
$L^{\otimes t}$ of the basic representation $L$
with the complex vector space $\Cc$ on which $1\in\Zz /p^{k+1}$ acts 
as multiplication by
$\zeta^t$ (where, we recall, $\zeta =\e^{2\pi\i /p^{k+1}}$). In particular, we can write $L^{\otimes t}=L^{\otimes t'}$
if $t$ and $t'$ are congruent modulo $p^{k+1}$.
There is an equivariant $t$th power map
$$
\pi_t : L\to L^{\otimes t}, \qquad z\mapsto z^t
$$
restricting to a map $S(L)\to S(L^{\otimes t})$ of unit circles.

We now need the following substantial input from stable homotopy 
theory.
\Prop{\label{meyercor}
{\rm (Stolz \cite{SS}), Meyer \cite{meyer}).}  
For $k=1$ and 
$d > 2$, there exists
a $\Zz /p^{k+1}$-map $S(p(d-2)L)\to S(dL^{\otimes p^k})$.
\qed
}
We shall manufacture maps for higher $k$ by using the $p$th power.
\Lem{\label{power}
Suppose that, for some integers $r,\, s\geq 1$
and $k\geq 1$,
there is a $\Zz /p^{k+1}$-map
$f : S(rL)\to S(sL^{\otimes p})$.
Let $\tilde L$ denote the basic $1$-dimensional complex representation
of $\Zz /p^{k+2}$.
Then there is a $\Zz /p^{k+2}$-map
$S(pr\tilde L)\to S(ps\tilde L^{\otimes p})$.
}
\begin{proof}
If $\Gamma$ is a finite group and $E$ is a Euclidean
$\Gamma$-module, then the $p$-fold
wreath product
$(\Gamma\times\cdots \times \Gamma)\semidirect \Zz /p$ acts on 
the join $S(E)*\cdots *S(E)$: 
an element $(g_1,\ldots ,g_p)\in \Gamma\times\cdots\times \Gamma$ maps
$[v_i,t_i]$, where $v_i\in S(E_i)$, $t_i\in [0,1]$, $\sum t_i=1$,
to $[g_iv_i,t_i]$ and $j\in\Zz /p$ sends
$[v_i,t_i]$ to $[v_{i+j},t_{i+j}]$.
We identify $S(E)*\cdots *S(E)$ with the sphere $S(\Sym^pE)$ on the
representation $\Sym^pE=E\oplus\cdots\oplus E$
(by mapping $[v_i,t_i]$ to $(\sqrt{t_i}\,v_i)$).
We shall make these constructions when $\Gamma$ is the
group $\Zz /p^{k+1}$.

\smallskip

Now the $p$th power of $f$
$$
f*\cdots * f : S(rL)*\cdots *S(rL) \to S(sL^{\otimes p})*
\cdots *S(sL^{\otimes p})
$$
is $(\Zz /p^{k+1}\times\cdots\times \Zz /p^{k+1})\semidirect 
\Zz /p$-equivariant.
We restrict from the action of the semidirect product to the action of
the cyclic subgroup $\Zz /p^{k+2}$
of order $p^{k+2}$ generated by $\tau =((1,0,\ldots ,0),1)
\in (\Zz /p^{k+1}\times\cdots\times \Zz /p^{k+1})\semidirect \Zz /p$.

Thus we obtain a $\Zz /p^{k+2}$-map 
$S(r\Sym^pL)\to S(s\Sym^p(L^{\otimes p}))$.
But $\Sym^p(L)$ is isomorphic, as $\Zz /p^{k+2}$-module, to
$\bigoplus_{i=0}^{p-1} \tilde L^{\otimes (1+ip^{k+1})}$
and $\Sym^p(L^{\otimes p})$ is isomorphic to
$\bigoplus_{i=0}^{p-1} \tilde L^{\otimes (p+ip^{k+1})}=
\bigoplus_{i=0}^{p-1} \tilde L^{\otimes p(1+ip^k)}$.

Composing $f*\cdots *f$ with the join of power maps
$$
\pi_{1+ip^{k+1}} :S(\tilde L)\to S(\tilde L^{\otimes (1+ip^{k+1})})
\text{\ and \ }
\pi_{1-ip^{k}}:
S(\tilde L^{\otimes p(1+ip^k)})\to S(\tilde L^{\otimes p}),
$$
since $(1+ip^k)(1-ip^k)\equiv 1$ (mod $p^{k+1}$), 
we produce the required map 
$S(pr\tilde L)\to S(ps\tilde L^{\otimes p})$.
\end{proof}
\Prop{\label{positive}
For $k\geq 1$, $0\leq  l \leq k$ and $d>2l$, 
there exists a $\Zz /p^{k+1}$-map
$$
S(p^k(d-2l)L) \to S(p^{k-l}dL^{\otimes p^{l}}).
$$
}
The statement is trivially true if $l=0$.
\begin{proof}
Starting from Proposition \ref{meyercor},
iteration using Lemma \ref{power} generates a $\Zz /p^{k+1}$-map
$S(p^k(d-2)L)\to S(p^{k-1}dL^{\otimes p})$ for all $k\geq 1$.
This establishes the result for $l=1$.

The argument is completed by induction on $l$. 

As in the proof of Lemma \ref{power} we write $\tilde L$ for the basic
complex representation of $\Zz /p^{k+2}$.
Suppose that we have a $\Zz /p^{k+1}$-map 
$S(p^k(d-2l)L) \to S(p^{k-l}dL^{\otimes p^{l}})$.
This may be regarded as a $\Zz /p^{k+2}$-map
$$
S(p^k(d-2l)\tilde L^{\otimes p}) \to 
S(p^{k-l}d\tilde L^{\otimes p^{l+1}}),
$$
which we can compose with a $\Zz /p^{k+2}$-map
$$
S(p^{k+1}((d-2l)-2)\tilde L) \to S(p^k(d-2l)\tilde L^{\otimes p})
$$
(as constructed as in the first paragraph)
to get a $\Zz /p^{k+2}$-map
$$
S(p^{k+1}(d-2(l+1))\tilde L) \to 
S(p^{k+1-(l+1)}d\tilde L^{\otimes p^{l+1}})
$$
if $d> 2(l+1)$. 
This completes the inductive step from $l$ to $l+1$.
\end{proof}
\begin{proof}[Proof of Theorem \ref{pos}]
If $k=0$, $V$ is isomorphic to $mL$ and we can take $c=0$.

Assume that $k\geq 1$. 
If $m_l\geq p^{k-l}(2l+1)$, write $m_l = n_lp^{k-l}+q_l$, where
$n_l\geq 1$ and $2lp^{k-l}\leq q_l <2(l+1)p^{k-l}$.
Otherwise, put $n_l=0$ and $q_l=m_l< p^{k-l}(2l+1)$.
Proposition \ref{positive} with $d=n_l+2l$ provides a map
$S(p^kn_lL) \to S(p^{k-l}(n_l+2l)L^{\otimes p^l})$.
Since $p^{k-l}(n_l+2l)\leq m_l$, we can compose with
the inclusion to get a $\Zz /p^{k+1}$-map
$$
S(p^kn_lL) \to  S(m_lL^{\otimes p^l}).
$$

Taking the join of such maps for $0\leq l\leq k$, we produce
a map $f_0:S(nL)\to S(V)$, where $n=\sum_l p^kn_l$.
Now
$$
n=\sum_{l=0}^k (p^lm_l -p^lq_l)\geq
(\sum_{l=0}^k p^lm_l) -\sum_{l=0}^k(2(l+1)p^k-p^l) =
(\sum_{l=0}^k p^lm_l) -c,
$$
where $c=p^k(k+2)(k+1)- (p^{k+1}-1)/(p-1)$.
\end{proof}
\begin{appendix}
\Sect{Cohomology and covering dimension}
In this appendix we collect for reference various standard facts related to the covering dimension of a compact Hausdorff space $Z$. Proofs may be found, for example, in \cite{kbu}.
\Lem{\label{bound}
Suppose that $Z$ is a closed subspace of a finite complex 
of dimension
$d$. Then $k^j(Z)=0$ and $H^j(Z;\,\Zz )=0$ for $j>d$.
\qed
}
\Lem{\label{coho}
Suppose that there is an integer $d$ such that $k^d(Z)$ is
non-zero and $k^j(Z)=0$ for $j>d$. 
Then the $\Zz$-cohomology group $H^d(Z;\,\Zz)$ is non-zero
and $H^j(Z;\,\Zz )=0$ for $j>d$.
\qed
}
\Lem{\label{coho2}
Suppose that there is an integer $d$ such that $(k\Zz /p^l)^d(Z)$ is
non-zero and $(k\Zz /p^l)^j(Z)=0$ for $j>d$. Then 
$H^j(Z;\,\Zz /p^l)=0$ for $j>d$ 
and $H^d(Z;\,\Zz /p^l)$ is non-zero.
Suppose further that there is an integer $c$ such that $H^j(Z;\,\Zz /p)=0$
for $j>c$.
Then $H^j(Z;\,\Zz /p)=0$ for $j>d$ and $H^d(Z;\,\Zz /p)$ is non-zero.
\qed
} 
Consider next a principal $\Zz /p$-bundle $\tilde Z \to Z$.
Let $\lambda$ be the associated complex line bundle
$\tilde Z\times_{\Zz /p}\Cc$ over $Z$, where the generator
$1\in\Zz /p$ acts as multiplication by $\omega =\e^{2\pi\i/p}$.
\Lem{\label{cover}
$\phantom{}$
\par\noindent
{\rm (i).}
Suppose that $H^j(\tilde Z;\,\Zz )=0$ for $j>d$.
Then $H^j(Z;\,\Zz )=0$ for $j>d$.

\par\noindent
{\rm (ii).} Suppose that $k^j(\tilde Z)=0$ for $j>d$.
Then $k^j(Z)=0$ for $j>d$.

There are similar results for cohomology and $K$-theory
with $\Zz /p^l$ coefficients for $l\geq 1$.
}
\begin{proof}
We prove (ii). The cover $\tilde Z$ is included fibrewise
in $S(\lambda)$ as 
$\tilde Z\times_{\Zz /p} \langle\omega\rangle$,
where $\langle\omega\rangle$ is the group of $p$th roots of
unity.
The exact sequence of the inclusion is
$$
\cdots\to
k^{j-1}(\tilde Z) \to k^j(S(\lambda )) \to k^j(\tilde Z) \to \cdots
$$
(because the quotient $S(\Cc )/\langle\omega\rangle$ is a wedge of
$p$ circles permuted by $\Zz /p$).
So $k^j(S(\lambda ))=0$ for $j>d+1$.

The Gysin exact sequence of the circle bundle $S(\lambda)$,
involving the $k$-theory Euler class $e_k(\lambda )$, is
$$
\cdots\to k^{j+1}(S(\lambda )) \to
k^{j}(Z)\Rarr{e_k(\lambda )}{} k^{j+2}(Z)
\to k^{j+2}(S(\lambda ))\to\cdots
$$
Hence, multiplication by $e_k(\lambda )$:
$k^j(Z)\to k^{j+2}(Z)$ is an isomorphism for $j>d$.
But $e_k(\lambda )$ is nilpotent.
\end{proof}
\Lem{\label{cdim}
Suppose that a compact Hausdorff space
$Z$ has finite covering dimension at most $d$. Then
$k^j(Z)=0$ and $H^j(Z;\,\Zz )=0$ for $j>d$.
\qed
}
\Prop{\label{orbit}
Let $\Gamma$ be a finite group and $Z$ a compact Hausdorff 
$\Gamma$-space.
Suppose that the compact Hausdorff orbit space $\bar Z=Z/\Gamma$ has covering
dimension at most $d$. Then $Z$ has covering dimension at most $d$.
\qed
}
\end{appendix}

\end{document}